\documentclass[12pt,twoside]{paper}
\usepackage{amsmath}
\usepackage{hyperref}
\usepackage{amsfonts}
\usepackage{amssymb}
\usepackage{amsthm}

\usepackage{graphicx}

\def\be{\begin{displaymath}}
\def\ee{\end{displaymath}}
\def\bee{\begin{equation}}
\def\eee{\end{equation}}

\def\Mas{Ma$\acute{\rm s}$lanka  }
\begin{document}

\thispagestyle{empty}
\centerline{}
\bigskip
\bigskip
\bigskip
\bigskip
\bigskip
\centerline{\Large\bf Evidence in favor of the}
\bigskip
\centerline{\Large\bf Baez-Duarte criterion for the Riemann
Hypothesis}
\bigskip

\begin{center}
{\large \sl Marek Wolf}\\*[5mm]

Institute of Theoretical Physics, University of Wroc{\l}aw\\
Pl.Maxa Borna 9, PL-50-204 Wroc{\l}aw, Poland, \href{mailto:mwolf@ift.uni.wroc.pl}{e-mail:mwolf@ift.uni.wroc.pl}\\
\bigskip

\end{center}

\bigskip
\bigskip

\bigskip\bigskip

\begin{center}
{\bf Abstract}\\
\bigskip
\bigskip
\begin{minipage}{12.8cm}
We presents results of the numerical experiments in favor of the Baez-Duarte criterion
for the Riemann Hypothesis.  We  give formulae allowing calculation of numerical values of the numbers $c_k$ appearing in this criterion for arbitrary large $k$. We present
plots of $c_k$ for $k \in (1, 10^9)$.
\end{minipage}

\end{center}
\bigskip\bigskip

{\bf 1. Introduction.}\\

In 1997 K. \Mas \cite{Maslanka1} proposed a new formula for the zeta Riemann
function valid on the whole complex plane $\mathbb{C}$ except a
point $s=1$:
\be
\zeta(s)=\frac{1}{1-s}\sum_{k=0}^\infty \frac{\Gamma(k+1-\frac{s}{2})}{\Gamma(1-\frac{s}{2})}
\frac{A_k}{k!}
\ee
where coefficients $A_k$ are given by
\be
A_k=\sum_{j=0}^k \binom{k}{j}(2j-1)\zeta(2j+2).
\ee

This formula was rigorously proved by L. Baez-Duarte
in 2003 \cite{Luis1}. In the subsequent preprint \cite{Luis2} the same author proved
the new criterion for the Riemann Hypothesis, the journal version of it appeared two
years later \cite{Luis3}.  The Riemann Hypothesis (RH) states that the nontrivial zeros $\rho$ of the function $\zeta(s)$ have the real part equal $\Re (\rho)=\frac{1}{2}$. Although Riemann
did not request it, today it is  often demanded additionally  that zeros on the critical line $\Re(s)=\frac{1}{2}$ should
be  simple. Baez-Duarte considered the sequence of numbers $c_k$ defined
by:
\bee
c_k=\sum_{j=0}^k {(-1)^j \binom{k}{j}\frac{1}{\zeta(2j+2)}}.
\label{ckmain}
\eee
He proved that RH is equivalent to the following rate of decreasing to zero of the above
sequence:
\bee
c_k={\mathcal{O}}(k^{-\frac{3}{4}+\epsilon})~~~~~~~~~{\rm for~ each~~}  \epsilon>0.
\label{criterion}
\eee
Furthermore, if $\epsilon$ can be put zero, i.e. if
$c_k={\mathcal{O}}(k^{-\frac{3}{4}})$, then the zeros of $\zeta(s)$ are simply. Baez-Duarte
also proved in \cite{Luis3} that it is not possible to replace $\frac{3}{4}$ by $\frac{3}{4}+
\epsilon$.

Neither in \cite{Luis3} nor in \cite{Maslanka3} it is explicitly
written whether the sequence $c_k$ starts from $k=0$ or $k=1$.
However in \cite{Luis3}  a few formulas contain
$k=0$, i.e.  summation starts from $c_0$. The point is that if we allow $k=0$, for which
$c_0=6/\pi^2$, then the inversion formula (see e.g. \cite{Knuth})
is fulfilled:

\bee
\frac{1}{\zeta(2k+2)}=\sum_{j=0}^k {(-1)^j \binom{k}{j} c_j}.
\label{inversion}
\eee
However I do not  see application of the above formula, except the possibility of checking
some of the statements made in \cite{Merlini}. Furthermore, if the Baez-Duarte sequence
$c_k$ starts from $k=0$ then the following identity holds:
\bee
\sum_{k=0}^\infty   \frac{c_k x^k}{k!} = e^x
\sum_{k=0}^\infty \frac{(-1)^k  x^k}{k!\zeta(2k+2)}.
\label{Cislo}
\eee
It is an application of the general formal identity:
\bee
\sum_{k=0}^\infty \left( \sum_{j=0}^k \binom{k}{j} a_j\right)
\frac{x^k}{k!} = e^x \sum_{k=0}^\infty \frac{a_k x^k}{k!},
\eee
where $a_k$ should not increase to fast with $k$ to ensure
convergence of series. \footnote{Indeed, collecting on the l.h.s.
terms multiplying $a_j$ we get: $ a_j \sum_{k=j}^\infty
\binom{k}{j} \frac{x^k}{k!} = a_j \sum_{k=j}^\infty \frac{k!}{j!
(j-k)!} \frac{x^k}{k!} =
a_j  \frac{x^j}{j!} \sum_{n=0}^\infty \frac{x^n}{n!} =  a_j
\frac{x^j}{j!} e^x$ and summing over $j$ gives r.h.s.}  Putting
here   $a_j=(-1)^j b_j$ gives the usual formula appearing in the
finite difference theory (see \cite{Flajolet} \S 1):
\bee
\sum_{k=0}^\infty \left( \sum_{j=0}^k (-1)^j \binom{k}{j}
b_j\right)  \frac{x^k}{k!} = e^x \sum_{k=0}^\infty \frac{(-1)^k
b_k x^k}{k!}.
\label{Cislo3}
\eee
The identity (\ref{Cislo}) can be used to
establish the connection with the Riesz criterion for RH (original
paper \cite{Riesz}, discussed in \cite{Luis3}, \cite{Luis4}).
Riesz has considered the function:
\begin{displaymath}
R(x) = \sum_{k=1}^\infty \frac{(-1)^{k+1}x^k}{(k-1)!\zeta(2k)}=
\sum_{k=0}^\infty \frac{(-1)^{k}x^{k+1}}{k!\zeta(2k+2)}.
\end{displaymath}
Unconditionally it can be proved that $R(x)=\mathcal{O}( x^{1/2+\epsilon})$, see
\cite{Titchmarsh}  \S 14.32. Riesz has proved that the Riemann Hypothesis is equivalent
 to slower increasing of the function $R(x)$:
\bee
RH \Leftrightarrow   R(x) = \mathcal{O}\left(
x^{1/4+\epsilon}\right).
\eee
But from (\ref{Cislo}) we get:
\bee
\sum_{k=0}^\infty   \frac{c_k x^k}{k!} = \frac {e^x}{x} R(x)
\label{Cislo2}
\eee
thus the generating function for $c_k$ can be
expressed by $R(x)$.
 In \cite{C-W} it is proved, that
for any real number $ \delta > -3/2 $ we have
\begin{equation}
R(x) = O(x^{\delta+1}) \Leftrightarrow  c_k = O(k^\delta).
\end{equation} \\
Proof is based on the relation $R(k)/k \approx c_k$.

\bigskip
\bigskip
{\bf 2. Computer experiments}
\bigskip

The criterion  (\ref{criterion}) seemed to be very well suited for the computer verification.
At the end of \cite{Luis2} Baez-Duarte wrote a sentence
``A test for the first $c_k$ up to $k = 1000$ shows a very pleasant smooth curve''.
However for larger values of $k$ the true behavior of the sequence turned out to be
more complicated: instead of monotonic tending to zero there appeared oscillations and
$c_k$ changed the sign at first for\footnote{in fact if the sequence $c_k$ starts from $k=0$ the
first sign change occurs for $c_0=6/\pi^2>0$ and $c_1=6/\pi^2-90/\pi^4=(6\pi^2-90)/\pi^4\approx(-30/\pi^4)<0$
(more precisely   $c_1=-0.3160113011...$)} $k=19320$: $c_{19319}=-1.7870567
\times 10^{-13}$
while $c_{19320}= 9.170232808\ldots
\times 10^{-12}$.   The next sign change is:  $c_{22526}=    2.2292905301\ldots
\times 10^{-13}$
but $c_{22527} =  -6.5057526
\times 10^{-12}$.

To my knowledge the first plot of $c_k$ for $k$ up to 95000 appeared in
the book \cite{Maslanka2} published in Polish. The same plot was reproduced in \cite{Luis3}.
Data used to make this plot consisted of $c_k$ calculated every 500-th $k$ --- it is very time
consuming to get $c_k$ directly from (\ref{ckmain}). Indeed, for large $j$ the values of
$\zeta(2j+2)$ very quickly become practically equal to 1, thus the summation of alternating
series gives wrong result when not performed with sufficient number of digits accuracy.
For example, the Table I below presents values of the partial sums for $c_{12000}$ recorded
every thousand summands (the calculation was performed with precision of 9000 digits).

Let us remark that the partial sums for $n$ and $12000 - n $ are of the same order.
The binomial coefficients become very large numbers in the middle and to get accurate
value of $c_k$ one needs a  lot of digits accuracy during the calculation. \Mas has
used {\tt  Mathematica } to perform these calculation. Over three years ago I
started to calculate
$c_k$ using  the free package PARI/GP \cite{Pari} developed especially for number theoretical
purposes and which allows practically arbitrary accuracy arithmetics both fixed-point as well
as floating-point. I started to calculate consecutive  $c_k$ for each $k$  with the help of the
following script in Pari:

\begin{verbatim}

\p 3500    /*     precision set to 3500 digits    */
allocatemem(250000000)
range=10000     /*     the largest subscript in c_k         */
denomin=vector(range);
for (n=1, range, denomin[n]=zeta(2*n));
default(format, "e22.20")
{
for  (k=1, range, c=sum(j=0,k,((-1)^j)*binomial(k,j)/denomin[j+1]);
write("c_k.dat",k," ", c))
}
\end{verbatim}

\vskip 0.4cm
\begin{center}
{\sf Table I} \\
\bigskip
\begin{tabular}{|c|c|} \hline
$ n $ & $ \sum_{j=0}^n {(-1)^j \binom{12000}{j}\frac{1}{\zeta(2j+2)}}$ \\ \hline
1000 &  8.6575528427959311728$\times 10^{1492}$ \\    \hline
2000  &    1.0610772171540382076$\times 10^{2346}$ \\    \hline
3000  &    2.6820721693716011525$\times 10^{2928}$ \\    \hline
4000  &    8.4511383022435967124$\times 10^{3314}$ \\    \hline
5000  &    1.8751018390471552047$\times 10^{3537}$ \\    \hline
6000  &    8.3417729099514988532$\times 10^{3609}$ \\    \hline
7000   &   1.3393584564622537177$\times 10^{3537}$ \\    \hline
8000  &    4.2255691511217983562$\times 10^{3314}$ \\    \hline
9000   &   8.9402405645720038417$\times 10^{2927}$ \\    \hline
10000   &   2.1221544343080764152$\times 10^{2345}$ \\    \hline
11000  &    7.8705025843599374298$\times 10^{1491}$ \\    \hline
12000  &    -1.6973092190852083930$\times 10^{-7}$ \\    \hline

\end{tabular}
\end{center}

\bigskip

The problem I have encountered during these calculations was that it seems to be not possible
to change accuracy of calculation during running the script (the command {\tt $\backslash$p 3500 }
above). Thus I had to change the precision by hand. It turned out that when the
precision was to small  produced values of $c_k$ were obviously wrong, something like
ten to the very large power. The rule learned from these
examples  for  precision  set to make calculations confident
was that the number of digits should be at least enough to distinguish between
1 and $1+\frac{1}{2^k}$ in the zeta appearing in (\ref{ckmain}), i.e. the precision set  to
calculate $c_k$ should be at least $\backslash {\tt p} =  k*\log_{10}(2)$. Table II presents
the real example I have met during calculations: when the precision was set to  60000 digits
the values of $c_k$ for k between 198000 and 200000 were  ($198000 \times \log_{10}(2)=
59603.93914,~~ 200000\times \log_{10}(2)=
60205.99913$):

\vskip 0.4cm
\begin{center}
{\sf Table II}\\
\bigskip
\begin{tabular}{|c|c|} \hline
$ k $ & $ c_k $ \\ \hline
198000 &    -8.1809420017968747912 $\times 10^{-9}$ \\ \hline
198500 & -8.1130397250007379108$\times 10^{-9} $\\ \hline
199000 &   -8.0431163120575296823$\times 10^{-9} $\\ \hline
199500 &  3.4122583912205353616 $\times 10^{49} $\\ \hline
200000 &  -1.9276608381598523688$\times 10^{200} $  \\ \hline

\end{tabular}
\end{center}

\Mas kindly send me values of $c_k$ from his calculations up to $k=95000$ with $k$ jumping
in intervals of 500, i.e. $k=500 l$. Autumn 2005 I have
started to continue this efforts on the cluster of 8 processors Xeon 2.8 GHz,  with
4 GB RAM per node of two processors\footnote{because I have used 32-bits version of PARI/GP
I was able to use $2^{31}$
bytes =2 GB  of RAM per process}, with the aim to reach
$k=200000$ also every 500-th value of $k$ using PARI/GP computer algebra system \cite{Pari}.
During last five months of computations between 4 and 6 processors I have
used to calculate $c_k$ in different intervals of $k$. When these calculations were running I
have learned of the paper \cite{Luis3} where explicit formulae for $c_k$ in terms
of zeros of $\zeta(s)$ were given. Quite recently there appeared the paper \cite{Maslanka3}
where the prescription to obtain $c_k$ very quickly were also given. In view of these
developments there is no need to continue very time consuming calculations based
on the formula (\ref{ckmain}). The only benefit of these calculation was the
possibility to compare $c_k$ obtained by means of formulae presented in \cite{Luis3}
and \cite{Maslanka3} against  those $c_k$ obtained from the generic formula (\ref{ckmain}).
It should be stressed that  calculations  based on (\ref{ckmain}) does not assume the validity
of Riemann Hypothesis in contrast to formulae presented by \Mas or below. Using these
formulae $c_k$ can be calculated very quickly for practically arbitrary $k$ ---
it is very time consuming to calculate $c_k$ without assuming RH.\\

\bigskip

{\bf 3. Explicit formulae.}\\

The formulae presented in \cite{Luis3}  and in \cite{Maslanka3} expressing $c_k$ directly
 in terms of the zeros
of $\zeta(s)$ are essentially the same, they differ in the manner they were derived.
\Mas has used the binomial transforms discussed  in \cite{Flajolet} while Baez-Duarte
is developing the whole machinery by himself. The formulae of these two authors can be
written as a sum of two parts: quickly decreasing with $k$ trend $\bar{c}_k$
 and oscillations $\tilde{c}_k$:

\be
c_k=\bar{c}_k +\tilde{c}_k
\ee
where:
\bee
\bar{c}_k=-\frac{1}{(2\pi)^2}\sum_{m=2}^\infty \frac{B(k+1,m)}{\Gamma(2m-1)}\frac{(-1)^m (2\pi)^{2m}}
{\zeta(2m-1)}
\label{cktrend}
\eee
and oscillating part:
\bee
\tilde{c}_k = \sum_{\rho} \frac{\Gamma(k+1)\Gamma(\frac{1+\rho}{2})}{\Gamma(k+1+\frac{1+\rho}{2})}
\frac{1}{\zeta'(1-\rho)} = \sum_{\rho} B\left(k+1,\frac{1+\rho}{2}\right) \frac{1}{\zeta'(1-\rho)}
\label{oscillations}
\eee
where it is assumed that zeros of $\zeta(s)$ are simple: $\zeta'(\rho)\neq 0$
and the sum is over all (i.e. on the positive as well as negative imaginary axis)
nontrivial zeros of the $\zeta(s)$, i.e. $\zeta(\rho)=0$
and $\Im \rho \neq 0$  and
\be
B(w,z)=\frac{\Gamma(w)\Gamma(z)}{\Gamma(w+z)}
\ee
is the Beta function. In fact Baez-Duarte is skipping the trend remarking  only that it
is of the order $o(1/k)$ (Remark 1.6 in \cite{Luis3}). Theoretically the formula for $\tilde{c}_k$
is valid in the limit of large $k$, but surprisingly the numbers produced from the above
formulae (\ref{cktrend}) and (\ref{oscillations}) are practically the same as obtained from
the generic formula (\ref{ckmain}) for all $k$, e.g. already for $k=2$ we get
$c_2=-0.25699711$ from
(\ref{ckmain}), while (\ref{cktrend}) and (\ref{oscillations}) give
$c_2= -0.256969863$  and accuracy increases with $k$. It suggests that the
integrals $J_k$ appearing in \cite{Luis3} in the proof of the Theorem 1.5 are decreasing
to zero rather fast with $k$.

First let us consider trend. It can be calculated directly from (\ref{cktrend}):
\bee
\bar{c}_k = -\frac{1}{(2\pi)^2} \sum_{m=2}^\infty \frac{1}
{(k+1)(k+2)\ldots(k+m)m(m+1)\ldots (2m-2)} \frac{(-1)^m (2\pi)^{2m}}{\zeta(2m-1)}
\label{exact}
\eee

\bigskip

\vskip 0.4cm
\begin{center}
{\sf Table III}\\
\bigskip
\begin{tabular}{|c||c|c|c|} \hline
$ k $ & $\bar{c}_k$ from eq.(\ref{exact})  & $\bar{c}_k$ from  eq.(\ref{trend_duze_k}) & $\bar{c}_k$ from  eq.(\ref{hugek})   \\ \hline
$ 1 $ & $ -2.60052406393\times 10^{-1 } $ &  $ -4.0752814729 \times 10^{3} $ &  $  -1.6421193331\times 10^{1} $ \\ \hline
$ 10 $ & $ -6.9069591105 \times 10^{-2}$  & $ -2.0455052855\times 10^{1} $ &  $  -1.6421193331\times 10^{-1} $   \\ \hline
$ 10^2 $  &  $   -1.4804264464\times 10^{-3} $ & $ -5.9943727867\times 10^{-3} $ &   $  -1.6421193331\times 10^{-3} $   \\ \hline
$ 10^{3} $ & $ -1.6248041420\times 10^{- 5} $ & $ -1.6824923096\times 10^{- 5} $ &  $ -1.6421193331\times 10^{-5} $ \\ \hline
 $ 10^{4} $ & $-1.6403755367\times 10^{- 7} $ &  $ -1.6418022398\times 10^{- 7} $ &  $ -1.6421193331\times 10^{-7 } $ \\ \hline
 $ 10^{5} $ & $ -1.6419448299\times 10^{- 9} $ & $  -1.6420436474\times 10^{- 9} $ & $ -1.6421193331\times 10^{-9} $ \\ \hline
 $ 10^{6} $ & $ -1.6421018816\times 10^{-11} $ & $ -1.6421113244\times 10^{-11} $ &  $ -1.6421193331\times 10^{-11} $ \\ \hline
 $ 10^{7} $ & $ -1.6421175880\times 10^{-13} $ & $ -1.6421185279\times 10^{-13} $  &  $ -1.6421193331\times 10^{-13} $ \\ \hline
 $ 10^{8} $ & $  -1.6421191586\times 10^{-15} $ & $ -1.6421192526\times 10^{-15} $  & $ -1.6421193331\times 10^{-15} $ \\ \hline
 $ 10^{9} $ & $  -1.6421193157\times 10^{-17}  $  &  $ -1.6421193251\times 10^{-17} $ & $  -1.6421193331\times 10^{-17} $ \\ \hline
\end{tabular}
\end{center}
\vskip 0.4cm

Using this formula I was able to produce every 500-th value of $\bar{c}_k$ for $k=500,1000,
\ldots 10^9$ performing calculations in Pari with 100 digits accuracy in about 4 hours.
For large $k$ I have used following asymptotic expansion of (\ref{exact}):
\bee
\bar{c}_k = -\frac{1}{(k+1)(k+2)}\left(\frac{(2\pi)^2}{2\zeta(3)}-
\frac{(2\pi)^4}{12(k+3)\zeta(5)}
+\frac{(2\pi)^6}{120(k+3)(k+4)\zeta(7)}\right).
\label{trend_duze_k}
\eee

It  can be further simplified to:
\bee
\bar{c}_k = -\frac{1}{k^2}\frac{(2\pi)^2}{2\zeta(3)}
\label{hugek}
\eee
The comparison of these formulae is given in Table III.

Now we consider the oscillating part $\tilde{c}_k$. Since PARI/GP does not have
built in $B(x,y)$ function, I had to use $\Gamma(z)$ functions instead.
Because of the fast growth of the $\Gamma(x)$ function even in PARI/GP it was not possible
to pursue with formula (\ref{oscillations}) for large $k$. Namely it crashes for $k=356000 $
because of overflow. But there is a following asymptotic formula (see e.g. \cite{Titchmarsh2},
\S 1.8.7):
\be
\frac{\Gamma(x)}{\Gamma(x+a)} \sim x^{-a}, ~~~~~~x \rightarrow \infty
\ee
thus we have
\bee
B(a,x)\sim x^{-a} \Gamma(a) ~~~~~~{\rm for}~~x~~{\rm large}.
\label{beta_duze}
\eee
Using it  for large $k$ and  assuming the Riemann Hypothesis: $\rho_l=\frac{1}{2}
+i\gamma_l$, $\bar{\rho}_l=\frac{1}{2} - i\gamma_l (= 1 - \rho_l)$ after collecting
together in pairs conjugate zeros  we get:
\bee
\tilde{c}_k=\frac{2}{(k+1)^{\frac{3}{4}}}\sum_{l=1}^\infty \alpha_l \cos\left(
\frac{1}{2}\gamma_l\log(k+1)\right)-
\beta_l\sin\left(\frac{1}{2}\gamma_l\log(k+1)\right),
\label{c_duze_k}
\eee
where I have denoted:
\begin{eqnarray}
\alpha_l=\Re \left(\frac{\Gamma(\frac{1+\rho_l}{2})}{\zeta^\prime(\bar{\rho_l})} \right), \\
\beta_l=\Im \left( \frac{\Gamma(\frac{1+\rho_l}{2})}{\zeta^\prime(\bar{\rho_l})} \right).
\end{eqnarray}
In (\ref{c_duze_k}) the decreasing of $c_k$ like $k^{-\frac{3}{4}}$ is obtained as an overall
amplitude of the ``waves'' composed of the cosines and sines with  the ''frequencies''
proportional to imaginary parts of the nontrivial zeros of $\zeta(s)$.
The coefficients $\alpha_l$ and $\beta_l$ decrease to zero very fast with $l$.
Namely using the Hadamard product for $\zeta(s)$:
\be
\zeta(s)=\frac{(2\pi)^se^{-(1+C/2)s}}{2(s-1)\Gamma(s/2+1)}\prod_\rho
\left(1-\frac{s}{\rho}\right)e^{\frac{s}{\rho}},
\ee
where $C=0.57721566490153286\ldots $ is the Euler constant,  the derivative of $\zeta(s)$ at
zeros can be computed. Taking into account miraculous simplifications, $\rho \rho_l +
\bar{\rho} \rho_l = \rho_l$ and the identity
\be
\Gamma(1+z)\Gamma(1-z)= \frac{\pi z}{\sin(\pi z)}
\ee
I have obtained that:
\bee
|\alpha_l| \varpropto e^{-\pi \gamma_l/4},~~~~~~|\beta_l|  \varpropto  e^{-\pi \gamma_l/4}.
\label{alpha_beta}
\eee
Because imaginary parts of zeros take large values it suffices to sum in (\ref{c_duze_k})
over a few first zeros. I have used 10 zeros and the table below gives  coefficients
$\alpha_l$ and $\beta_l$ and comparison with (\ref{alpha_beta}).

\vskip 0.4cm
\begin{center}
{\sf Table IV}\\
\bigskip
\begin{tabular}{|c||c|c|c|c|} \hline
$ l $ &  $ \alpha_l $  &   $\beta_l $ &  $e^{-\pi \gamma_l/4} $   \\ \hline
$ 1 $ & 2.029173866 $  \times 10^{-5} $ & $ -3.315924256  \times 10^{-5} $       &  $ 1.50914 \times 10^{-5} $  \\   \hline
$ 2 $ & -3.333265938 $  \times 10^{-8} $ &$ -1.298336420   \times 10^{-7}$       & $ 6.75315 \times 10^{-8} $  \\   \hline
$ 3 $ & 2.886139424 $  \times 10^{-9} $ & $ -4.153918097   \times 10^{-9}$       & $ 2.94404 \times 10^{-9} $  \\   \hline
$ 4 $ & 4.813880001 $  \times 10^{-11} $ & $  -6.332017430   \times 10^{-11} $   &  $ 4.19039 \times 10^{-11} $  \\   \hline
$ 5 $ & 7.546769513 $  \times 10^{-12} $ & $ 7.526891498  \times 10^{-12}$       & $ 5.83506 \times 10^{-12} $  \\   \hline
$ 6 $ & 6.162524600 $  \times 10^{-14} $ & $  1.942118979  \times 10^{-13}$      & $ 1.51209 \times 10^{-13} $  \\   \hline
$ 7 $ & -1.578482027 $  \times 10^{-14} $ & $ 1.184829593   \times 10^{-14}$     & $ 1.10374 \times 10^{-14} $  \\   \hline
$ 8 $ & -1.07138189 2$  \times 10^{-15} $ & $  -2.209146437  \times 10^{-15}$    & $ 1.66491 \times 10^{-15} $  \\   \hline
$ 9 $ & 9.328038737 $  \times 10^{-19} $ & $ -7.472197226  \times 10^{-17}$      & $ 4.22403 \times 10^{-17} $  \\   \hline
$ 10 $ & 1.747829093 $  \times 10^{-17} $ & $ 1.122667624 \times 10^{-17}$       & $ 1.05303 \times 10^{-17} $ \\   \hline
\end{tabular}
\end{center}
\vskip 0.4cm

However already calculations with the first zero $\gamma_1=14.13472514173469\ldots$
give numbers which differ
much less than 1\%  (see $l=1$ and $l=2$ in the above table) from  those calculated with larger
number of terms in (\ref{c_duze_k}) as well as  with $c_k$  for $k<200000$
 obtained directly from (\ref{ckmain}) without assuming RH.  
The plots of $c_k$ for $k$ up to $10^9$ obtained from these formulae are given in the Fig.1
and Fig.2. In Fig.2 there is  logarithmic  $k$-axis and thus the plot has a constant
``wavelength'', not depending on $k$ like on the Fig.1.  The envelope is given by:
\bee
y=\frac{2A}{(k+1)^{3/4}}~,~~~~~~A=0.777506276445256\times 10^{-5}
\label{envelope}
\eee
and was obtained in the following way: 
First I have maintained in (\ref{c_duze_k}) only the first zero $\rho_1=\frac{1}{2}+\gamma_1$:
\bee
\tilde{c}_k=\frac{2}{(k+1)^{\frac{3}{4}}}\left( \alpha_1 \cos\left(
\frac{1}{2}\gamma_1\log(k+1)\right)-
\beta_1\sin\left(\frac{1}{2}\gamma_1\log(k+1)\right)\right).
\eee
Next I made use of the identity:
\bee
a \cos(\theta) - b \sin(\theta) = A \sin(\phi-\theta)~~~{\rm where}~~~A=\sqrt{a^2+b^2},
~~~\phi=\arctan\left(\frac{a}{b}\right)
\eee
to obtain:
\bee
\tilde{c}_k=\frac{2}{(k+1)^{\frac{3}{4}}}\sqrt{\alpha_1^2+\beta_1^2}\sin\left(\phi -\frac{1}{2}\gamma_1\log(k+1)\right).
\label{sinus}
\eee
from which (\ref{envelope}) follows and numerical value of $A$ is obtained from $\alpha_1$ and
$\beta_1$  in Table IV. Here$\phi=-0.54916(=56.497^\circ)$ Let us remark that this value of $A$ agrees very well with amplitude
reported by Beltraminelli and Merlini \cite{Merlini2}. It is interesting to note that lhs of the  above  formula is valid not only for integer $k$ but also for real $k$, thus using the approximation $c_k\approx R(k)/k$ derived in \cite{C-W} we can write for large $x$:
\bee
R(x)=2x^{1/4}\left( \alpha_1 \cos\left(\frac{1}{2}\gamma_1\log(x)\right)-
\beta_1\sin\left(\frac{1}{2}\gamma_1\log(x)\right)\right)
\eee

There is another way of checking accuracy of the above equation (\ref{sinus}). Namely assuming that (\ref{sinus})
 is true and denoting  by $k^{'}$ and $k^{''}$ two consecutive zeros of $c_{k'}=0, c_{k''}=0$
we get for $\gamma_1$

\bee
\gamma_1=\frac{2\pi}{\log((k^{''}+1)/(k^{'}+1))}
\label{gamma}
\eee
To make sense, in the latter approach an independent  of (\ref{sinus}) and relatively fast method of calculating
$c_k$ is needed. In fact in \cite{Luis3} Baez-Duarte gives among others following
formula being the transformation of (\ref{ckmain}):
\bee
c_k=\sum_{n=1}^{\infty} \frac{\mu(n)}{n^2}\left(1-\frac{1}{n^2}\right)^k.
\label{c k moebius}
\eee

Using this expression I have searched numerically for sign changes of $c_k$ up
to $k=10^{9}$ and Table V presents $\gamma_1$
calculated from the above formula (\ref{gamma})  together with consecutive zeros extrapolated from
integer values using the linear approximation (let us remind that $\sin(x)$ has derivative
$\pm 1$ at zeros!). The correct value is $\gamma_1= 14.13472514173469379045725\ldots$.

In Fig.2 the plot of  $ k^{\frac{3}{4}} c_k$ is presented. The Baez-Duarte
criterion requires this  ``wave'' to be contained in the strip of parallel lines for all $k$.
The violation of the RH would manifest as an increase of the amplitude of the combination
$ k^{\frac{3}{4}} c_k$ for large $k$. This point is elaborated in more detail
by \Mas in \cite{Maslanka3}. Here I will make some further comments on this issue.
First it should be   remarked that the r.h.s. of (\ref{c_duze_k}) consists of products of
three terms: the first depending only on $k$ (the overall factor $k^{\frac{3}{4}}$), the
second depending only on imaginary parts of nontrivial zeros of $\zeta$ (the coefficients
$\alpha_l$ and $\beta_l$) and third ingredients depending both on $k$ and $l$ (the
trigonometric functions). Assume there are some zeros of $\zeta$ off critical line.
 We can split the sum over zeros $\rho$ in
(\ref{oscillations}) in two parts: one over zeros on critical line and second over zeros off
critical line.  This second sum should  violate the overall term $k^{-3/4}$ present in
the first sum.  Let $\gamma_i^{(o)}$ denote  the imaginary parts of the zeros
lying off critical line (``o'' stands for ``off''). It is not clear whether asymptotic similar
to (\ref{alpha_beta}) will be valid for zeros off critical line, but it seems to be
reasonable to assume that it should not differ significantly from (\ref{alpha_beta}).
Then the  contribution to $c_k$ of such zeros off critical line should contain a factor of
the order $e^{-\gamma_l^{(o)}}$. Because value of the
imaginary part $\gamma_l^{(o)}$ of the hypothetical zero off critical line should be extremely
large, perhaps even as large as $10^{100}$ can be expected, the combined contribution to $c_k$
coming from the second sum seems to be extremely small, thus to see violation of the
Baez-Duarte criterion the values of $k$ should be larger than famous Skewes number and
look something like $10^{{10}^{{.}^{{.}^{.}}}}$. Such a big index $k$ should cause that
the first term in (\ref{c_duze_k}) overcame the
smallness of the second term depending only on $\gamma_l^{(o)}$.

\vskip 0.4cm
\begin{center}
{\sf Table V}\\
\bigskip
\begin{tabular}{|c||c|c|} \hline
  & extrapolated zero $k^{'}$  & $\gamma_1$  from eq. (\ref{gamma}) \\ \hline
1 &       19319.0191151 &   0.66394362155812867 \\ \hline
2 &       22526.0331312 &  40.91301517418643030 \\ \hline
3 &       41868.9707418 &  10.13657613078084811 \\ \hline
4 &       60094.3311655 &  17.38744852979003230 \\ \hline
5 &       98378.2514809 &  12.74744116303923005 \\ \hline
6 &      149320.1629630 &  15.05785832239914317 \\ \hline
7 &      236817.0977574 &  13.62376898016027559 \\ \hline
8 &      366000.4553802 &  14.43265811387491946 \\ \hline
9 &      573460.6753253 &  13.99205468779138406 \\ \hline
10 &      891841.3774543 &  14.22828231751450665 \\ \hline
11 &     1393469.2943691 &  14.07955314904129825 \\ \hline
12 &     2173554.0482344 &  14.13327453400666383 \\ \hline
13 &     3387835.5708183 &  14.15682029348025265 \\ \hline
14 &     5283842.8916393 &  14.13659885512530467 \\ \hline
15 &     8247263.1465316 &  14.11229699162269810 \\ \hline
16 &    12864372.4128156 &  14.13284964864925287 \\ \hline
17 &    20052822.4883780 &  14.15424390183719722 \\ \hline
18 &    31271608.9210745 &  14.14046804681939058 \\ \hline
19 &    48805962.1935733 &  14.11501826908551253 \\ \hline
20 &    76145459.4891092 &  14.12608973663948854 \\ \hline
21 &   118689214.0783221 &  14.15568694594979885 \\ \hline
22 &   185076299.0519995 &  14.14304477168657435 \\ \hline
23 &   288851950.0033488 &  14.11488330063446422 \\ \hline
24 &   450663149.7663923 &  14.12567920473796440 \\ \hline
25 &   702517003.8056590 &  14.15292988626078933 \\ \hline
\end{tabular}
\end{center}
\bigskip

The plot in Fig.2 is a  perfect sine of one wavelength 
thus it gives visual  justification of the above statement that $\tilde{c}$ is determined
in fact by the first zero $\gamma_1$. The same phenomenon was mentioned by \Mas in
\cite{Maslanka3}.

\bigskip
\bigskip

{\bf 4. Final remarks}

\bigskip

The formula (\ref{c_duze_k}) together with a few coefficients $\alpha_l$ and $\beta_l$ taken
from Table IV allows to compute values of $c_k$ for  arbitrary large $k$. Other criteria for
RH, like the value of the  de Bruijn-Newman constant \cite{Odlyzko2}, are vulnerable to
the Lehmer pairs of zeros of $\zeta(s)$. It is hard to see the reason for violation of
the inequality $|c_k| < {\rm const}~ k^{-\frac{3}{4}}$. I have checked, that at the first Lehmer pair
$\rho_{6709}=0.5+ 7005.06286617 i$ and $\rho_{6710}=0.5 + 7005.1005646 i$
the derivative has value $\zeta'(0.5+ 7005.06286617 i) = 3.2229849698 + 0.74179951875 i$
and similar value for second zero,
thus there is no chance to get values of $c_k$ violating (\ref{criterion}) in this way.
 It seems to
be an open problem how to  connect the value of the largest $k$ for which
 $|c_k| < {\rm const} k^{-\frac{3}{4}}$ to the number of zeros lying on the critical line.
Let us mention that for the Li's criterion \cite{Li} which states that if the  numbers
\bee
\lambda_n=\frac{1}{(n-1)!}\frac{d^n}{ds^n}(s^{n-1}\log \xi(s)) |_{s=1} 
\eee
fulfill $\lambda_n>0$ for each $n$ then RH is true it is known that if the first $n$ Li's constants
$\lambda_n$ are positive
then every zero $\rho$ of $\zeta(s)$ with $|\Im \rho |<\sqrt{n}$ lies on the critical
line $\Re \rho= \frac{1}{2}$ \cite{AMS}.

After a few months of computer experiments with $c_k$ I believe Baez-Duarte sequence
is one the most important  and mysterious sequences in the whole mathematics.

\bigskip
\bigskip

{\bf Acknowledgement} I thank Dr J.Cis{\l}o for discussions and bringing to my attention
identity (\ref{Cislo3}), Prof. L. Baez-Duarte and Prof. K. \Mas for e-mail exchange. The zeros
of $\zeta(s)$ I have used in computations were taken from Prof. A. Odlyzko's web page
\href{http://www.dtc.umn.edu/~odlyzko/zeta\_tables/index.html}
{http://www.dtc.umn.edu/$\sim$odlyzko/zeta\_tables/index.html}


\begin{center}
  \begin{figure}
    \includegraphics[width=1\columnwidth, clip=true]{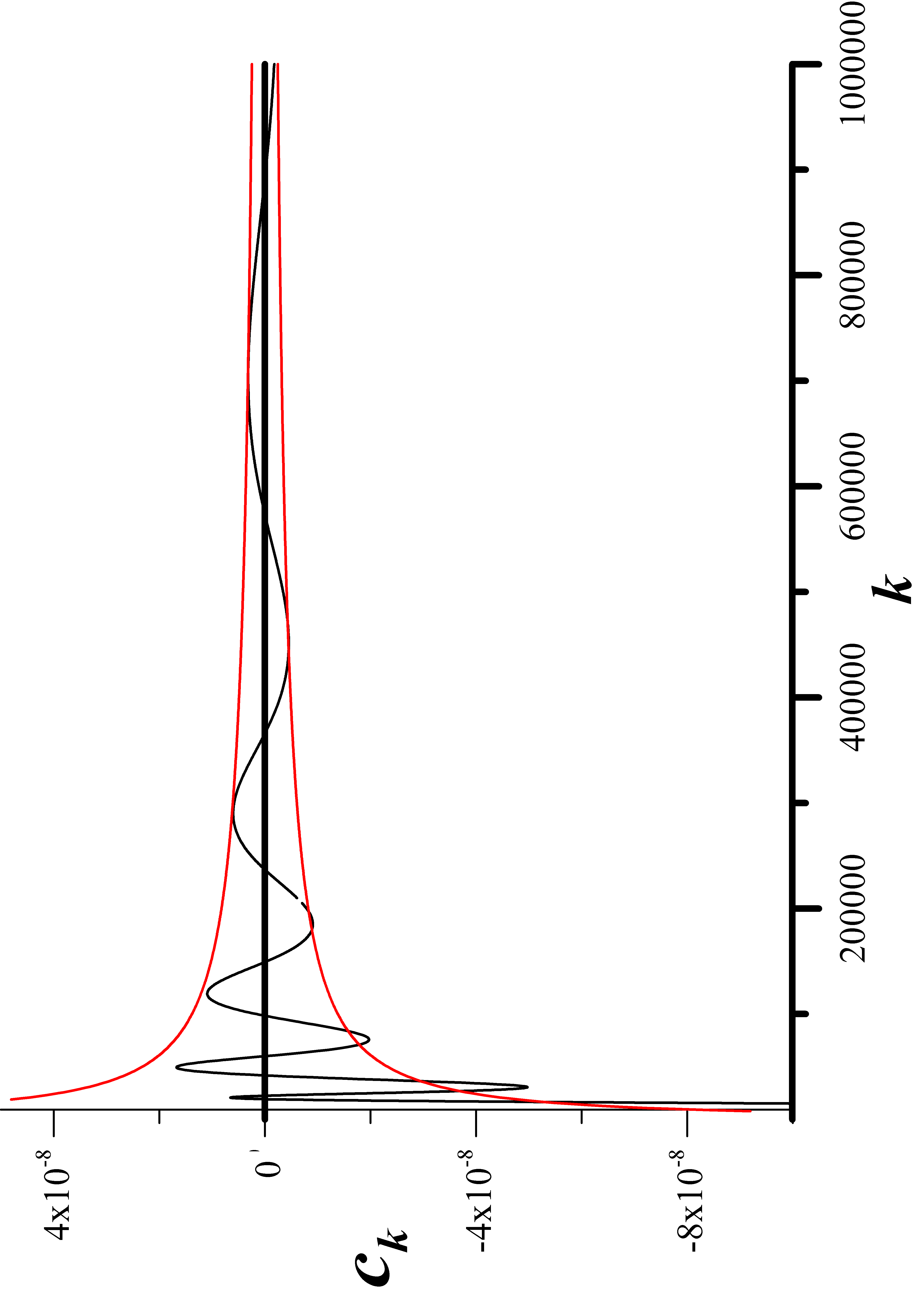}
    \caption{{The plot of $c_k$ for $k\in(1, 10^6)$. At $k=200000$ a small gap is visible to distinguish between $c_k$ calculated from  generic formula  (\ref{ckmain}) and from
    explicit formulas presented in Sec. 3.}
    \label{fig:1}}
  \end{figure}
\end{center}


\begin{center}
  \begin{figure}
    \includegraphics[width=1\columnwidth, scale=1.2, clip=true]{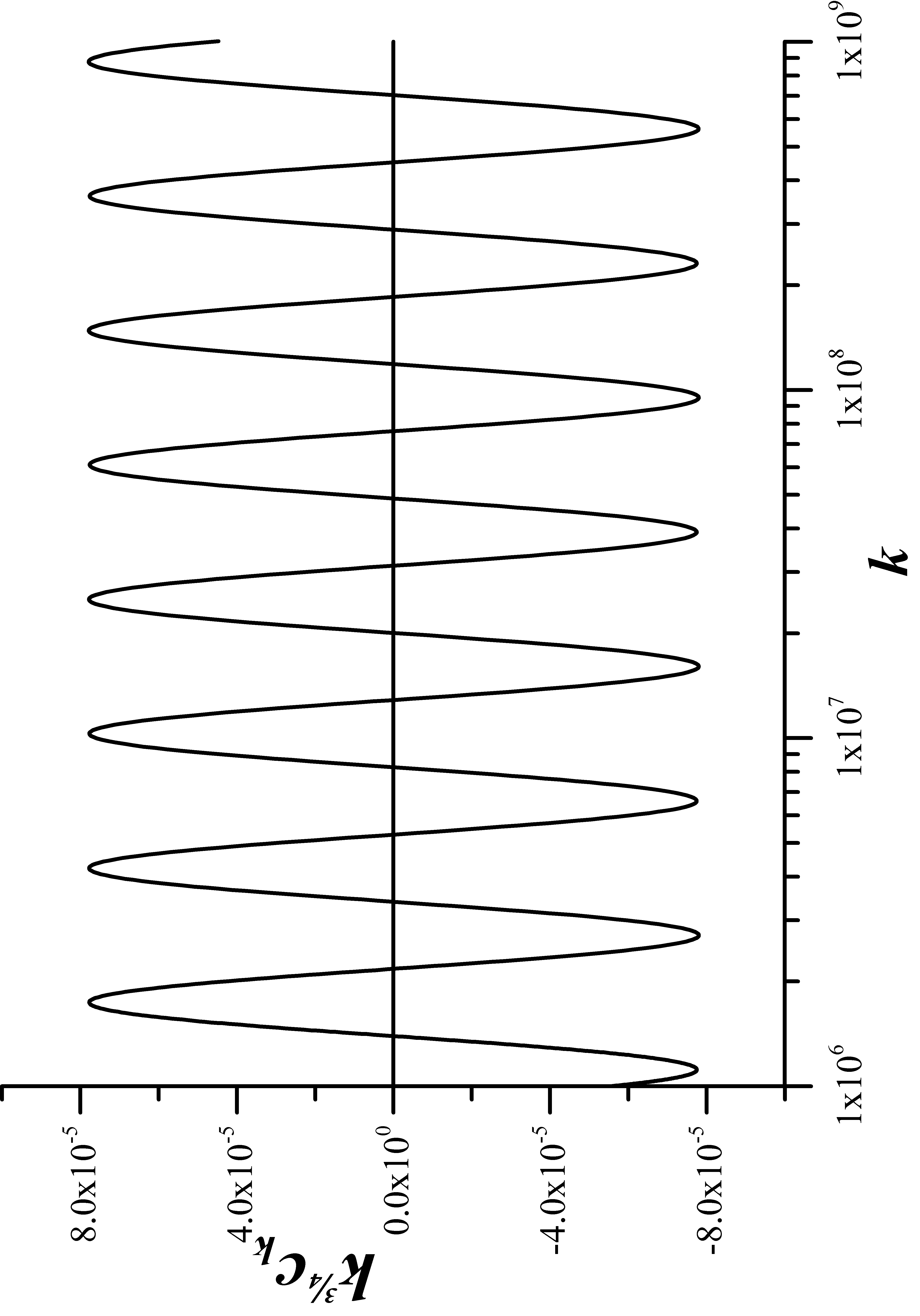}
    \caption{{ The plot of $k^{3/4}c_k$ for $k\in(10^6,10^9)$.}
    \label{fig:1}}
  \end{figure}
\end{center}

\end{document}